\documentclass[journal]{IEEEtran}

\ifCLASSINFOpdf
\else
\fi

\usepackage{array}
\usepackage{amssymb,amsmath}
\usepackage{graphicx}
\usepackage{amsthm}
\usepackage{epstopdf}
\allowdisplaybreaks[4]
\usepackage{amsmath}
\usepackage{flushend}
\usepackage{color}
\usepackage{makecell}
\usepackage{float}
\emergencystretch=\maxdimen
\newtheorem{theorem}{Theorem}
\newtheorem{definition}{Definition}

\newtheorem{lemma}{Lemma}
\newtheorem{assumption}{Assumption}
\newtheorem{remark}{Remark}

\newcommand{\col}{$\upshape{col}$}

\newcommand{\diag}{$\upshape{diag}$}

\begin{document}
%
\title{A distributed generalized Nash equilibrium seeking algorithm based on extremum seeking control}
%
%
%

\author{\IEEEauthorblockN{Feng~Xiao, {\it Member, IEEE}, Xin~Cai, and Bo~Wei}

\thanks{
This work was supported by the National Natural Science Foundation of China (NSFC, Grant  Nos. 61873074, 61903140).
(Corresponding author: Feng~Xiao.)}

\thanks{F. Xiao and X. Cai are with the State Key Laboratory of Alternate Electrical Power System with Renewable Energy Sources and with the School of Control and Computer Engineering, North China Electric Power University, Beijing 102206, China (emails: fengxiao@ncepu.edu.cn; caixin\_xd@126.com). X. Cai is also with the School of Electrical Engineering, Xinjiang University, Urumqi 830047, China.}

\thanks{B. Wei is with the School of Control and Computer Engineering,
North China Electric Power University,
Beijing 102206, China (email: bowei@ncepu.edu.cn).}}

\maketitle

\begin{abstract}
In this paper, a distributed non-model based seeking algorithm which combines the extremum seeking control (ESC) jointly with learning algorithms is proposed to seek a generalized Nash equilibrium (GNE) for a class of noncooperative games with coupled equality constraint. The strategy of each agent is restricted by both the coupled inter-agent constraint and local inequality constraints. Thanks to the ESC, it is unnecessary to know the specific expressions of agents' cost functions and local constraints and to know the strategies of other agents for the implementation of the proposed GNE seeking algorithm. To deal with the coupled constraints, only the Lagrange multiplier is transmitted among agents with some prior information about the coupled constraints. Moreover, a diminishing dither signal is designed in the seeking algorithm to remove undesirable steady-state oscillations. The non-local convergence of the designed seeking algorithm is analyzed via the singular perturbation theory, averaging analysis and Lyapunov stability theory. Numerical examples are given to verify the effectiveness of our proposed method.
\end{abstract}


%
\IEEEpeerreviewmaketitle

\section{Introduction}
In network scenarios, distributed control of multi-agent systems based on game theory have attracted an increasing interest in various engineering applications, such as mobile sensor networks \cite{stankovic.2012,Zhang.2018,Ma.2015,Wang.2016,Yu.2021,Ma.2019}, communication networks \cite{Wang.2014,Yi.2019}, smart grids \cite{Gharesifard.2016,Li.2021,Qin.2021}, energy markets \cite{Wang.2021,Mu.2012,Guo.2016}, and so on. In such a case, the control objective of multi-agent systems corresponds to the Nash equilibrium (NE) of games. Thus, many researchers have paid much attention to seek the Nash equilibrium, especially for generalized Nash equilibrium (GNE) of games with coupled constraints.

Most of the existing GNE seeking algorithms were designed based on the primal-dual theory. The designed algorithms mainly include two parts. One is to update the strategies of agents by the gradient of a Lagrange function which is composed of the cost function, Lagrange multipliers and constraints. The other is the coordination of Lagrange multipliers which are used to drive strategies to satisfy coupled constraints. According to various application scenarios, existing distributed algorithms to seek GNE can be roughly categorised into two classes: discrete-time algorithms \cite{Yi.2019,Liu.2020,Yi.2020,Yi.2019a,Lu.2021} and continuous-time algorithms (\cite{Lu.2019,Liang.2017,Zeng.2019,Ye.2017} and references therein), which we focus on in this paper.  Based on consensus protocols, which were used to estimate some necessary global information, distributed algorithms were designed to seek GNE of $N$-person noncooperative games for agents with smooth cost functions on undirected communication graphs \cite{Lu.2019,Ye.2017}. Since communication networks play an important role in the design and implementation of distributed algorithms, the GNE seeking problem was extended from fixed undirected graphs to directed and time-varying graphs \cite{Deng.2019b,Liang.2017}. In addition, distributed algorithms were proposed to deal with nonsmooth aggregative games and nonsmooth multi-cluster games, respectively \cite{Deng.2021,Zeng.2019}. To remove the requirement of some information of communication graphs, adaptive gains were introduced to distributed GNE seeking algorithms \cite{Bianchi.2019,Li.2020}.

It is worth noting that the above-mentioned work require the exact information of the model, i.e., the specific expressed cost functions and constraints, and the strategies of other players. However, in some practical engineering systems, there exist unknown cost functions, constraints and other parameters in games, such as unknown price functions in energy markets \cite{Ye.2016}, the power function without explicit aerodynamic interactions in the wind turbine \cite{Marden.2013}, unknown traffic demands/constraints in the network routing \cite{Marden.2009}, and unknown position information and/or environmental conditions in mobile sensor networks \cite{stankovic.2012}. To deal with such situations, an adaptive and real-time method, called extremum seeking control (ESC) \cite{Tan.2006,Tan.2008,Tan.2009,Tan.2013}, was introduced in distributed seeking algorithms \cite{stankovic.2012,Frihauf.2012,Ye.2020,Zahedi.2019}. In this paper, we explore how to steer agents to update their strategies adaptively to GNE under the conditions of unknown cost functions and local constraints. For $N$-person noncooperative static games, the classical ESC method was used to solve the NE \cite{stankovic.2012,Frihauf.2012,Ye.2020}. To eliminate steady-state oscillations, an extremum seeking method with fast convergence was designed to seek the NE \cite{Zahedi.2019}. The comparison between the existing algorithms and our proposed method is shown in Table \ref{Table 1}. To the best of our knowledge, the GNE seeking for games with coupled constrains remains to be solved. Due to local inequality constraints and coupled equality constraint, existing distributed algorithms (see \cite{Yi.2019,Yi.2020,Lu.2019} and references therein) are invalid to seek NE of games in the case that only the values of cost functions is available. Moreover, classical ESC based distributed algorithms designed in \cite{stankovic.2012,Frihauf.2012,Zahedi.2019,Ye.2020} also fail to seek GNE of games with coupled constraints. Therefore, these observations motivate us to work out an extremum seeking method for noncooperative games with coupled constraints.

\begin{table} [H]
\caption{Comparison with the existing NE seeking algorithms}
\scalebox{0.87}{
\label{Table 1}
\centering
\begin{tabular}{c|cccc}
\hline
\makecell[c]{NE seeking algorithms\\
based on ESC} & \makecell[c]{steady-state\\ oscillations} & \makecell[c]{local\\ constraints} & \makecell[c]{coupled\\ constraints} & \makecell[c]{non-local\\ convergence}\\
\hline
\cite{stankovic.2012,Frihauf.2012} & $\checkmark$ & $\times$ & $\times$ & $\times$ \\{}
\cite{Ye.2020} & $\checkmark$ & $\times$ & $\times$ & $\checkmark$ \\{}
\cite{Zahedi.2019} & $\times$ & $\times$ & $\times$ & $\times$ \\{}
our proposed method & $\times$ & $\checkmark$ & $\checkmark$ & $\checkmark$ \\
\hline
\end{tabular}}
\end{table}

In this paper, a distributed extremum-seeking based GNE seeking algorithm is designed. Compared with existing work, the main contributions of this paper are summarized as follows.

1) A distributed GNE seeking algorithm is proposed based on the ESC with only access to the values of cost functions and local constraints. The local inequality constraints are dealt with an exact penalty method which attributes to the values of the local constraints used in the ESC. The coupled equality constraint is coped with the Lagrange multiplier which is regulated by a consensus estimator in a distributed manner. Due to the essential feature of multiple time scales of the ESC, it is technically challenging to analyze the designed extremum seeking algorithm combined with the regulation of Lagrange multipliers.

2) To remove the undesirable steady-state oscillation existing in the classical ESC, inspired by \cite{Tan.2009,Wang.2016}, the amplitude of dither signal is regulated adaptively in the designed algorithm. In contrast to the local convergence obtained in \cite{Zahedi.2019}, the designed algorithm ensures non-local convergence.

The rest of this paper is organized as follows. In Section \uppercase\expandafter{\romannumeral2}, the considered problem is formulated. In Section \uppercase\expandafter{\romannumeral3}, a GNE seeking algorithm based on extremum seeking is designed. In Section \uppercase\expandafter{\romannumeral4}, special cases of general quadratic functions and stubborn players are analyzed. In Section \uppercase\expandafter{\romannumeral5}, simulation examples are presented. The conclusions are stated in Section \uppercase\expandafter{\romannumeral6}.

\section{Problem Formulation}
\subsection{Notations}
In this paper, $\mathbb{R}$ and $\mathbb{R}_{\geq0}$ denote the set of real numbers and non-negative real numbers, respectively. $\mathbb{R}^n$ is the set of $n$-dimensional real vectors. $\mathbb{R}^n_{\geq0}$ denotes the set of $n$-dimensional real vector with non-negative elements. $\mathbb{R}^{n\times m}$ denotes the set of $n \times m$ real matrices. Given a vector $x\in \mathbb{R}^n$, $\|x\|$ is the Euclidean norm. $A^T$ and $\|A\|$ are the transpose and the spectral norm of matrix $A\in \mathbb{R}^{n\times n}$, respectively. $\lambda_i(A)$ is the $i$th eigenvalue of matrix $A$, and it is expressed simply as $\lambda_i$. $\diag\{h_1,\ldots,h_n\}$ is a diagonal matrix with diagonal elements $h_1,\ldots,h_n$. $\boldsymbol{1}_n$ and $\boldsymbol{0}_n$ are $n$-dimensional column vectors with all elements being ones and zeros, respectively. $I_n$ denotes the $n\times n$ identity matrix. $\boldsymbol{0}$ is a zero matrix with an appropriate dimension. The gradient of $g(x)$ with respect to $x$ is denoted by $\nabla g(x)$. Additionally, $[u]^{+}=\max\{0,u\}$ for $u\in\mathbb{R}$. A continuous function $\sigma:[0,a)\rightarrow[0,\infty)$ is a class $\mathcal{K}$ function if it is strictly increasing and $\sigma(0)=0$. A continuous function $\beta(r,t):[0,a)\times[0,\infty)\rightarrow [0,\infty)$ is a class $\mathcal{KL}$ function if, for each fixed $t$, $\beta(r,t)$ is a class $\mathcal{K}$ function with respect to $r$ and, for each fixed $r$, $\beta(r,t)$ is decreasing with respect to s and $\beta(r,t)\rightarrow 0$ as $t\rightarrow \infty$.

\subsection{Problem Formulation}

In this paper, we consider an $N$-person noncooperative game with coupled constraints, which is denoted by $G=(\mathcal{I},\Omega,J,C)$. Each player is indexed by the set $\mathcal{I}=\{1,\ldots,N\}$. The strategy space $\Omega$ of the game is denoted by $\Omega=\Omega_1\times \cdots \times \Omega_N$, where $\Omega_i\subset\mathbb{R}$ is the player~$i$'s strategy set characterized by local inequality constraints $g_{ij}(x_i)\leq 0$ $(j=1,\ldots, m_i)$ with $g_{ij}: \mathbb{R}\rightarrow \mathbb{R}$. The cost profile $J$ of the game is constituted by the cost functions $J_i(x_i,x_{-i}): \mathbb{R}\times\mathbb{R}^{N-1}\rightarrow\mathbb{R}$ of player $i$,  $i\in\mathcal{I}$, where $x_i$ is the strategy of player $i$, $x_{-i}$ $=$ $\col(x_1,\ldots,x_{i-1},x_{i+1},\ldots,x_N)$ denotes the strategy vector of players except player $i$. $C$ denotes the coupled constraint characterized by $\sum_{i=1}^Nx_i=\sum_{i=1}^Nd_i$. Given the strategies of other players, player $i$ chooses his strategy $x_i$ by solving the following optimization problem.
\begin{equation} \label{op}
\begin{array}{l}
\min_{x_i} J_i(x_i,x_{-i})\\
s.t.\ \   g_{ij}(x_i)\leq 0, \ \ j=1,\ldots,m_i\\
\ \ \ \sum_{i=1}^Nx_i=\sum_{i=1}^Nd_i.
\end{array}
\end{equation}

Suppose that players can communicate with each other via an undirected and connected graph $\mathcal{G}=(\mathcal{I},\mathcal{E})$, and the game admits a GNE $x^*$ which is a strategy profile denoted by $(x_1^*,\ldots,x_N^*)$. For an undirected and connected graph $\mathcal{G}$, the sets of nodes and edges are denoted by $\mathcal{I}$ and $\mathcal{E}$, respectively. And there exists a path between any pair of distinct vertices. The Laplacian matrix $L$ of an undirected and connected graph $\mathcal{G}$ has the following properties \cite[Lemma 6.4; Theorem 6.6]{Bullo.2020}: (1) $\lambda_1(L)=0<\lambda_2\leq \cdots \leq \lambda_N$; (2) $\boldsymbol{1}_N^TL=\boldsymbol{0}_N$ and $L \boldsymbol{1}_N=\boldsymbol{0}_N$.The objective of this paper is to propose a GNE seeking algorithm for the game $G$, given that players can only measure the values of their own cost functions and local constraints.

\begin{definition}[Generalized Nash
Equilibrium, GNE {\cite[Definition 3.7]{Basar.1999}}, \cite{Facchinei.2007}] A strategy profile $x^*=[x_1^*,\ldots,x_N^*]^T$ is said to be a generalized Nash equilibrium of the game $G$ if
\begin{equation*}
J_i(x_i^*,x_{-i}^*)\leq J_i(x_i,x_{-i}^*),  \forall x_i\in\Omega_i\cap C, i\in\mathcal{I}.
\end{equation*}
\end{definition}

\begin{figure*}[!t]
\normalsize
\begin{equation} \label{ghj}
\nabla_{x_i}\hat{J}_i(x_i,x_{-i})=
\begin{cases}
\nabla_{x_i}J_i(x_i,x_{-i}), &  g_{ij}(x_i)\leq0, \forall j\in \{1,\ldots,m_i\}\\
\nabla_{x_i}J_i(x_i,x_{-i})+p_i\sum_{j=1}^{m_i}\nabla g_{ij}(x_i), &  g_{ij}(x_i)>0, \forall j\in \{1,\ldots,m_i\}
\end{cases}.
\end{equation}
\hrulefill
\vspace*{4pt}
\end{figure*}

\begin{assumption} \label{as1}
The cost function $J_i(x_i,x_{-i})$  and local constraints $g_{ij}(x_i)$  $(j=1,\ldots,m_i)$ are convex in $x_i$. And they are sufficiently smooth for all $i\in\mathcal{I}$. Generalized Nash equilibrium to the problem \eqref{op} exists and satisfies the Slater's condition.
\end{assumption}

\begin{assumption} \label{as2}
For any different $x$, $y\in\mathbb{R}^N$,
\begin{equation*}
(x-y)^T(F(x)-F(y))\geq m\|x-y\|^2,
\end{equation*}
where $F(x)=[\nabla_{x_1} J_1(x_1,x_{-1}),\ldots,\nabla_{x_N} J_N(x_N,x_{-N})]^T$ and $\nabla_{x_i} J_i(x_i,x_{-i})$ denotes the gradient of $J_i(x_i,x_{-i})$ with respect to $x_i$.
\end{assumption}

Assumption \ref{as1} ensures that the problem \eqref{op} is convex. Assumption \ref{as2} characterizes a unique GNE in the game $G$. Note that Assumptions \ref{as1}-\ref{as2} about cost functions and local constraints indicate that the extremum seeking scheme is only applicable to certain types of cost functions, although the scheme does not require knowing the explicit expression of cost functions.

The exact penalty method for inequality constraints is utilized to convert the problem \eqref{op} of agent $i$ to the following auxiliary problem.

\begin{equation} \label{eop}
\begin{array}{l}
\min_{x_i} \hat{J}_i(x_i,x_{-i})\\
s.t.\ \  \sum_{i=1}^Nx_i=\sum_{i=1}^Nd_i,
\end{array}
\end{equation}
where $\hat{J}_i(x_i,x_{-i})=J_i(x_i,x_{-i})+p_i\sum_{j=1}^{m_i}[g_{ij}(x_i)]^{+}$ is called an auxiliary cost function, and the penalty coefficient $p_i>0$. The gradient of $\hat{J}_i(x_i,x_{-i})$ with respect to $x_i$ is expressed in Eq.~\eqref{ghj}.

If local constraints are differential convex functions in $x_i$, the gradient of $g_{ij}(x_i)$ is monotone for all $j\in\{1,\ldots,m_i\}$, that is, $(x_i-x_i')(\nabla g_{ij}(x_i)-\nabla g_{ij}(x_i'))>0$ for any different $x_i$, $x_i'\in\mathbb{R}$ and for all $j\in\{1,\ldots,m_i\}$. It is easy to conclude that under Assumption \ref{as1}, the auxiliary cost function $\hat{J}(x_i,x_{-i})$ satisfies Assumption \ref{as2} if $J_i(x_i,x_{-i})$ satisfies it.

Under Assumption \ref{as1}, for fixed $x_{-i}$, the problems \eqref{op} and \eqref{eop} are convex optimization problems. The relationship between problems \eqref{op} and \eqref{eop} is given as follows.

\begin{lemma}[{\cite[Proposition 2.2]{Cherukuri.2015}},\cite{Boyd.2004}] \label{lemma1} Under Assumption \ref{as1}, the problem \eqref{op} is convex and strong duality holds. Then, any solution of the problem \eqref{eop} is also a solution of \eqref{op}, if $p_i>\max\{\eta_{i1},\ldots,\eta_{im_{i}}\}$, where $\eta_{ij}\in\mathbb{R}_{\geq 0}$ is the Lagrange multiplier corresponding to local inequality constraint $g_{ij}\leq 0$ $(j=1,\ldots,m_i)$ in the problem \eqref{op}.
\end{lemma}

Note that for all $i\in \mathcal{I}$, the Lagrange multiplier $\eta_{ij}$ $(j=1,\ldots,m_i)$ in problem \eqref{op} exists and satisfies KKT conditions. The penalty coefficient $p_i$ can be selected sufficiently large without the need of computing $\eta_{ij}$ $(j=1,\ldots,m_i)$.

\begin{lemma}[{\cite[Theorem 3.1]{Facchinei.2007}}] \label{lemma2}
Under Assumptions \ref{as1}-\ref{as2}, for problem \eqref{eop}, there exists a unique GNE $x^*=[x_1^*,\ldots,x_N^*]^T$ and a common Lagrange multiplier $\bar{\mu}\in\mathbb{R}$ such that
\begin{align*}
\nabla_{x_i} J_i(x_i^*,x_{-i}^*)+\bar{\mu}&=0,\\
\sum_{i=1}^N x_i^*=\sum_{i=1}^N d_i,
\end{align*}
for all $i\in\mathcal{I}$.
\end{lemma}

\begin{remark}
For the sake of simplicity, we only consider the case where the strategy of each player is scalar. The proposed method can be applied to the $m$-dimensional vector case by introducing $m$ virtual players whose cost functions are the same as the one of the corresponding real player \cite{Frihauf.2012}.
\end{remark}

\section{Main results}
\begin{figure}[H]
  \centering
  \includegraphics[width=8cm]{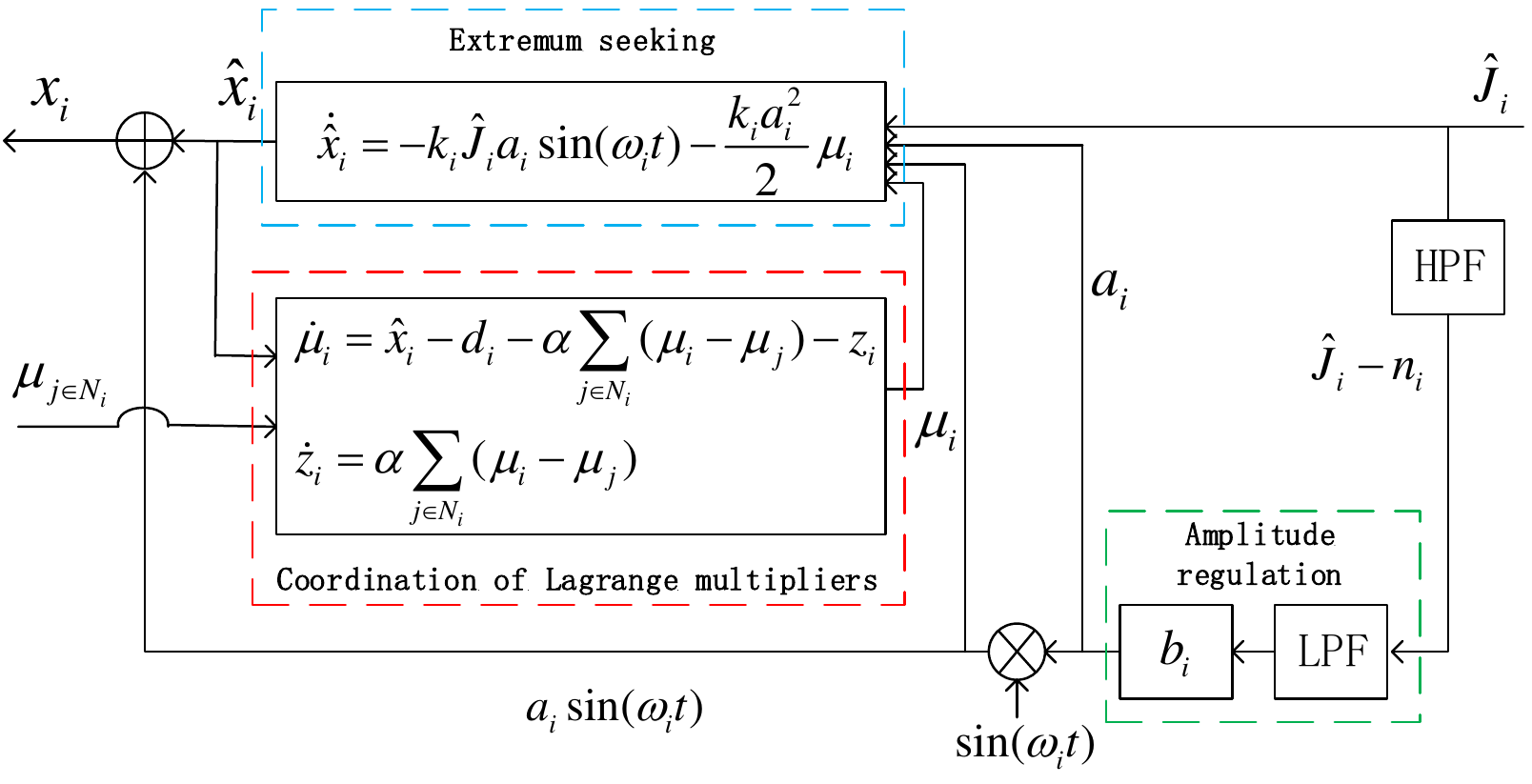}
  \caption{The proposed GNE seeking algorithm}
  \label{fig3}
\end{figure}

In this section, a novel GNE seeking algorithm is designed. The structure of the designed algorithm is shown in Fig.~\ref{fig3}. The designed algorithm mainly includes three parts, that is, the extremum seeking for the estimation of the auxiliary cost function's gradient, the coordination of the Lagrange multipliers and the regulation of the dither signal's amplitude.

The GNE seeking algorithm for agent $i$ is proposed as follows.
\begin{align} \label{sur}
\begin{split}
x_i&=\hat{x}_i+a_i\sin(w_it),\\
\dot{\hat{x}}_i&=-k_i\hat{J}_ia_i\sin(w_it)-\frac{k_ia_i^2}{2}\mu_i,\\
\dot{\mu}_i&=\hat{x}_i-d_i-\alpha\sum_{j\in\mathcal{N}_i}(\mu_i-\mu_j)-z_i,\\
\dot{z}_i&=\alpha\sum_{j\in\mathcal{N}_i}(\mu_i-\mu_j),\\
\dot{a}_i&=-\delta w_{li}a_i+\delta b_iw_{li}(\hat{J}_i-n_i),\\
\dot{n}_i&=-w_{hi}n_i+w_{hi}\hat{J}_i,
\end{split}
\end{align}
where $a_i$ and $w_i=w\bar{w}_i$ are the amplitude and frequency of the dither signal, respectively. $w$ is a positive real constant and $\bar{w}_i$ is some fixed positive constant. $\hat{J}_i$ is the auxiliary cost in problem \eqref{eop} and its value can be obtained by the measurements of cost $J_i$ and local constraints $g_i$. $\mu_i$ is the Lagrange multiplier corresponding to the coupled constraint. $z_i$ is an auxiliary variable with initial state $\sum_{i=1}^Nz_i(0)=0$ or $z_i(0)=0,\forall i\in\mathcal{I}$. $n_i$ is the low frequency component of $\hat{J}_i$. $k_i$, $b_i$ and $\alpha$ are positive constants. $\delta$ is a small positive constant. $w_{li}$ and $w_{hi}$ are cut off frequencies of low-pass and high-pass filters, respectively. As shown in Fig \ref{fig3}, HPF and LPF are the high-pass and low-pass filters, respectively. Let $\underline{K}=\min\{k_1a_1^2,\ldots,k_Na_N^2\}$. We have the following conclusion.

\begin{theorem} \label{theorem1}
Let $D$ be a compact set that contains the origin. Suppose that Assumptions \ref{as1}-\ref{as2} hold, the parameters $\underline{K}$ and $\alpha$ satisfy that $\underline{K}>\frac{1}{2\phi m-1}$ and $\alpha>\frac{4+\phi^2m}{2\lambda_2}$ with $\phi>\frac{1}{2m}$, respectively, and the frequencies of dither signals satisfy that $w_i\neq w_j$, $w_i\neq 2w_j$, $w_i\neq w_j+w_k$, $2w_i\neq w_j+w_k$ and $w_i\neq 2w_j+w_k$ for any distinct $i,j,k\in\mathcal{I}$. All the players follow the seeking algorithm \eqref{sur} to update their strategies. Then, there exist $w^*>0$ and $\delta^*>0$ such that for any $w>w^*$, $\delta\in(0,\delta^*)$,
\begin{align*}
\|\Xi(t)\|\leq \beta(\|\Xi(t_0)\|,(t-t_0))+O(1/w), \ \ \forall t\geq t_0 \geq 0
\end{align*}
with $\Xi(t_0)\in D$, $\|\Xi(t_0)\|\leq\Delta$ and $a_i(t_0)>0, i\in\mathcal{I}$, where $\beta$ is a $\mathcal{KL}$ function, $\Xi(t)=\col(x-x^*,\mu-\mu^*,z-z^*,a-a^*,n-n^*)$, and $\Delta$ depends on $\max\{a_1(t_0),\ldots,a_N(t_0)\}$.
\end{theorem}

%
%

Proof: Let $\sigma=\delta \tau$ with $\tau=wt$. System \eqref{sur} can be rewritten in $\sigma$-time scale.
\begin{align} \label{sps}
\begin{split}
\delta \frac{d\hat{x}_i}{d\sigma}&=\frac{1}{w}(-k_i\hat{J}_ia_i\sin(w'_i\sigma)-\frac{k_ia_i^2}{2}\mu_i),\\
\delta \frac{d\mu_i}{d\sigma}&=\frac{1}{w}(\hat{x}_i-d_i-\alpha\sum_{j\in\mathcal{N}_i}(\mu_i-\mu_j)-z_i),\\
\delta \frac{dz_i}{d\sigma}&=\frac{1}{w}(\alpha\sum_{j\in\mathcal{N}_i}(\mu_i-\mu_j)),\\
\frac{da_i}{d\sigma}&=\frac{1}{w}(-w_{li}a_i+b_iw_{li}(\hat{J}_i-n_i)),\\
\delta \frac{dn_i}{d\sigma}&=\frac{1}{w}(-w_{hi}n_i+w_{hi}\hat{J}_i),
\end{split}
\end{align}
where $w'_i=\bar{w}_i/\delta$.
It is easy to know that \eqref{sps} is a singular perturbation system with the parameter $\delta$. We set $\delta=0$ such that $(\hat{x}_i, \mu_i, z_i, n_i)$ is at the quasi-steady state. The obtained reduced system is given by
\begin{align} \label{rd}
\frac{da_i}{d\sigma}=-\frac{w_{li}}{w}a_i.
\end{align}
Its equilibrium point is denoted by $a_i^*=0$. It is concluded that $a_i$ exponentially converges to zero with $a_i(0)>0$.

For system \eqref{sur}, let $\delta=0$ freeze $a_i$ to yield the boundary layer system, given by
\begin{align} \label{bl}
\begin{split}
\frac{d\hat{x}_i}{d\tau}&=\frac{1}{w}(-k_i\hat{J}_ia_{i}\sin(\bar{w}_i\tau)-\frac{k_ia_{i}^2}{2}\mu_i),\\
\frac{d\mu_i}{d\tau}&=\frac{1}{w}(\hat{x}_i-d_i-\alpha\sum_{j\in\mathcal{N}_i}(\mu_i-\mu_j)-z_i),\\
\frac{dz_i}{d\tau}&=\frac{1}{w}(\alpha\sum_{j\in\mathcal{N}_i}(\mu_i-\mu_j)),\\
\frac{dn_i}{d\tau}&=\frac{1}{w}(-w_{hi}n_i+w_{hi}\hat{J}_i).
\end{split}
\end{align}

Since the sinusoidal perturbation is periodic, averaging theory can be utilized to analyze the system \eqref{bl}. The Taylor expansion of $\hat{J}_i(x)$ \footnote{With slight abuse of notation, $\hat{J}_i(x_i,x_{-i})$ is denoted by $\hat{J}_i(x)$.} at $\hat{x}$ is
\begin{align*}
\begin{split}
\hat{J}_i(x)&=\hat{J}_i(\hat{x})+\sum_{j=1}^N\frac{\partial \hat{J}_i(\hat{x})}{\partial \hat{x}_j}(a_j\sin(\bar{w}_j\tau))\\ &\ \ \ +\frac{1}{2}\!\sum_{k=1}^N\!\sum_{j=1}^N\frac{\partial^2 \hat{J}_i(\hat{x})}{\partial \hat{x}_k\partial \hat{x}_j}(a_k\sin(\bar{w}_k\tau))(a_j\sin(\bar{w}_j\tau))\\
&\ \ \ + a_i^3R_{i},
\end{split}
\end{align*}
where $R_{i}$ denotes the residual term of the Taylor expansion.


Then, if $w_i \neq w_j$, $w_i\neq 2w_j$, $w_i \neq w_j+w_k$, $2w_i \neq w_j+w_k$, and $w_i \neq 2w_j+w_k$, the average of boundary layer system \eqref{bl} is computed as follows.
\begin{align} \label{ave}
\begin{split}
\frac{d\hat{x}_i^a}{d\tau}&=\frac{1}{w}(-k_i\frac{a_{i}^2}{2}\frac{\partial \hat{J}_i(\hat{x}^a)}{\partial \hat{x}_i^a}-\frac{k_ia_{i}^2}{2}\mu_i),\\
\frac{d\mu_i}{d\tau}&=\frac{1}{w}(\hat{x}_i^a-d_i-\alpha\sum_{j\in\mathcal{N}_i}(\mu_i-\mu_j)-z_i),\\
\frac{dz_i}{d\tau}&=\frac{1}{w}(\alpha\sum_{j\in\mathcal{N}_i}(\mu_i-\mu_j)),\\
\frac{dn_i}{d\tau}&=\frac{1}{w}(-w_{hi}n_i+w_{hi}\hat{J}_i),
\end{split}
\end{align}
where $\hat{x}^a$ $=$ $[\hat{x}_1^a,\ldots,\hat{x}_N^a]^T$.

Let $\mu=$ $[\mu_1,\ldots,\mu_N]^T$, $z=[z_1,\ldots,z_N]^T$,
$n=$ $[n_1,\ldots,n_N]^T$, $d$ $=$ $[d_1,\ldots,d_N]^T$, $\hat{J}$ $=$ $[\hat{J}_1,\ldots,\hat{J}_N]^T$, $K$ $=$ $\diag \{k_1a_{1}^2,\ldots,k_Na_{N}^2\}$, $w_h$ $=$ $\diag\{w_{h1},\ldots,$ $w_{hN}\}$, and $F(\hat{x}^a)=[\frac{\partial \hat{J}_1(\hat{x}^a)}{\partial \hat{x}_1^a},\ldots,\frac{\partial \hat{J}_N(\hat{x}^a)}{\partial \hat{x}_N^a}]^T$. The compact form of system \eqref{ave} is given by
\begin{align}\label{cf}
\begin{split}
\frac{d\hat{x}^a}{d\tau}&=\frac{1}{w}(-\frac{K}{2}(F(\hat{x}^a)+\mu)),\\
\frac{d\mu}{d\tau}&=\frac{1}{w}(\hat{x}^a-d-\alpha L\mu-z),\\
\frac{dz}{d\tau}&=\frac{\alpha}{w}L\mu,\\
\frac{dn}{d\tau}&=\frac{1}{w}(-w_h(n-\hat{J})).
\end{split}
\end{align}
The equilibrium of system \eqref{cf} is denoted by $(x^*,\mu^*,z^*,n^*)$, where $x^*$ is the GNE of the game $G$ by Lemma \ref{lemma2}, $\mu^*$$=$$\bar{\mu}\boldsymbol{1}_N$ for some $\bar{\mu}\in \mathbb{R}$ and $n^*=J^*$.

Consider the following Lyapunov function $V=w\phi(\hat{x}^a-x^*)^T(\hat{x}^a-x^*)+\frac{w}{2}(\mu-\mu^*)^T(\mu-\mu^*)+\frac{w}{2}(\mu-\mu^*+z-z^*)^T(\mu-\mu^*+z-z^*)+\frac{w}{2}(n-n^*)^Tw_h^{-1}(n-n^*)$
with $\phi>0$. The derivative of $V$ along the trajectories of \eqref{cf} is given by
\allowdisplaybreaks[4]
\begin{align} \label{dly}
\begin{split}
\frac{dV}{d\tau}&=-\phi(\hat{x}^a-x^*)^TK(F(\hat{x}^a)-F(x^*)+\mu-\mu^*)\\
&\ \ \ +(\mu-\mu^*)^T(\hat{x}^a-\alpha L\mu-z-x^*+\alpha L\mu^*+z^*)\\
&\ \ \ +(\mu-\mu^*+z-z^*)^T(\hat{x}^a-z-x^*+z^*)\\
&\ \ \ -(n-n^*)^T(n-n^*)\\
&\leq-\phi\underline{K}(\hat{x}^a-x^*)^T(F(\hat{x}^a)-F(x^*))\\
&\ \ \ -\phi\underline{K}(\hat{x}^a-x^*)^T(\mu-\mu^*)+2(\mu-\mu^*)^T(\hat{x}^a-x^*)\\
&\ \ \ -\alpha(\mu-\mu^*)^TL(\mu-\mu^*)-2(\mu-\mu^*)^T(z-z^*)\\
&\ \ \ +(z-z^*)^T(\hat{x}^a-x^*)-\|z-z^*\|^2-\|n-n^*\|^2,
\end{split}
\end{align}
where $\underline{K}=\min\{k_1a_1^2,\ldots,k_Na_N^2\}$.

Each term in \eqref{dly} is analyzed as follows. According to Assumption 2, it yields that
\begin{equation} \label{t1}
-\phi\underline{K}(\hat{x}^a-x^*)^T(F(\hat{x}^a)-F(x^*))\leq -\phi\underline{K}m\|\hat{x}^a-x^*\|^2.
\end{equation}
Recall that the communication graph $\mathcal{G}$ is undirected and connected. It is easy to get
\begin{equation} \label{t2}
-\alpha(\mu-\mu^*)^TL(\mu-\mu^*)\leq -\alpha\lambda_2\|\mu-\mu^*\|^2.
\end{equation}
Using Young's inequality, we have that
\begin{equation} \label{t3}
-\phi\underline{K}(\hat{x}^a-x^*)^T(\mu-\mu^*)\leq \frac{\underline{K}}{2}\|\hat{x}^a-x^*\|^2+\frac{\phi^2\underline{K}}{2}\|\mu-\mu^*\|^2.
\end{equation}
Besides,
\begin{equation} \label{t4}
\begin{split}
&2(\mu\!-\!\mu^*)^T(\hat{x}^a\!-\!x^*)\!-\!2(\mu\!-\!\mu^*)^T(z\!-\!z^*)\!+\!(z\!-\!z^*)^T(\hat{x}^a\!-\!x^*)\\
&= -\|-\sqrt{2}(\mu-\mu^*)-\frac{\sqrt{2}}{2}(z-z^*)+\frac{\sqrt{2}}{2}(\hat{x}^a-x^*)\|^2\\
&\ \ \ +2\|\mu-\mu^*\|^2+\frac{1}{2}\|z-z^*\|^2+\frac{1}{2}\|\hat{x}^a-x^*\|^2.
\end{split}
\end{equation}
Substituting \eqref{t1} - \eqref{t4} into \eqref{dly} leads to
\begin{align*}
\begin{split}
\frac{dV}{d\tau}
&\leq -(\phi\underline{K}m-
\frac{\underline{K}}{2}-\frac{1}{2})\|\hat{x}^a-x^*\|^2-\frac{1}{2}\|z-z^*\|^2\\
&\ \ \ -(\alpha\lambda_2-\frac{\phi^2\underline{K}}{2}-2)\|\mu-\mu^*\|^2-\|n-\hat{J}^*\|^2.
\end{split}
\end{align*}
If $\phi>\frac{1}{2m}$, $\underline{K}>\frac{1}{2\phi m-1}$ and $\alpha>\frac{4+\phi^2m}{2\lambda_2}$, $\dot{V}<0$. The average of $(\hat{x}_i,\mu_i,z_i,n_i)$-subsystem can asymptotically converge to $(x_i^*,\bar{\mu},z_i^*,n_i^*)$ for all $i\in\mathcal{I}$.

Next, we will show that how the behavior of average system \eqref{ave} approximates the behavior of the boundary layer system \eqref{bl}. Define
$h_i(\tau,\!\hat{x}_i)\!=\!-k_i\big(\!\hat{J}_i(\hat{x})\sin(\!w_i\tau\!)\!+\!\frac{\partial \!\hat{J}_i(\hat{x})}{\partial \!\hat{x}_i}a_{i}\!\sin^2\!(w_i\tau)\!\big)\!+\!\frac{k_i}{2}\!\frac{\partial \!\hat{J}_i(\hat{x})}{\partial \! \hat{x}_i}$
for $i={1,\ldots,N}$. It is obvious that $h_i(\tau,\hat{x}_i),\forall i\in \mathcal{I}$, is periodic in $\tau$ and has zero mean for $\hat{x}_i$, which is in the compact set $D$. Define $H_i(\tau,\hat{x}_i)=\int_0^{\tau} h_i(s,\hat{x}_i)ds$, $\forall i\in\mathcal{I}$. Then, $H_i(\tau,\hat{x}_i)$ is periodic in $\tau$ and bounded for $\hat{x}_i$ which is in the compact set $D$. In addition,
\begin{equation}
\begin{aligned}
\frac{\partial H_i(\tau,\hat{x}_i)}{\partial \tau}&=h_i(\tau,\hat{x}_i),\\
\frac{\partial H_i(\tau,\hat{x}_i)}{\partial \hat{x}_i}&=\int_0^{\tau} \frac{\partial h_i(s, \hat{x}_i)}{\partial \hat{x}_i}ds,
\end{aligned}
\end{equation}
 are periodic in $\tau$ and bounded for $\hat{x}_i$. Define a new coordinate $y_i$ according to
 \begin{align} \label{nc}
\hat{x}_i=y_i+a_{i}H_i(\tau,y_i).
 \end{align}
Differentiating both sides of \eqref{nc} with respect to $\tau$ yields that
\begin{equation*}
\frac{d\hat{x}_i}{d\tau}=\frac{dy_i}{d\tau}+a_{i}(\frac{\partial H_i(\tau,y_i)}{\partial \tau}+\frac{\partial H_i(\tau,y_i)}{\partial y_i}\frac{dy_i}{d\tau}).
\end{equation*}
It derives that
\begin{equation*}
\begin{split}
(1+a_{i}\frac{\partial H_i(\tau,y_i)}{\partial y_i})\frac{dy_i}{d\tau}&= \frac{d\hat{x}_i}{d\tau}-a_{i}\frac{\partial H_i(\tau,y_i)}{\partial \tau}\\
&=-\frac{k_ia_{i}^2}{2}(\frac{\partial \hat{J}_i(y)}{\partial y_i}+\mu_i)+k_ia_{i}q_i,
\end{split}
\end{equation*}
where $y=[y_1,\ldots,y_N]^T$, and $q_i$ $=$ $-(\hat{J}_i(\hat{x})\sin(w_i\tau)$ $+$ $\frac{\partial \hat{J}_i(\hat{x})}{\partial \hat{x}_i}a_{i}\sin^2(w_i\tau))+\hat{J}_i(y)\sin(w_i\tau)+\frac{\partial \hat{J}_i(y)}{\partial y_i}a_{i}\sin^2(w_i\tau)+O(a_{i}^2)$. By the mean value theorem, it can be derived that $q_i=O(a_{i})+O(a_{i}^2)$ for $y$ in any compact set, $i\in\mathcal{I}$. Since $1+a_{i}\frac{\partial H_i}{\partial y_i}\neq 0$, $(1+a_{i}\frac{\partial H_i}{\partial y_i})^{-1}=1+O(a_{i})$ for sufficiently small $a_{i}$.

Thus, $\frac{dy_i}{d\tau}=-\frac{k_ia_{i}^2}{2}(\frac{\partial \hat{J}_i(y)}{\partial y_i}+\mu_i)+O(a_{i}^2)+O(a_{i}^3)$. With the new coordination $y_i$, the system \eqref{bl} can be rewritten as follows.
\begin{align} \label{ncs}
\begin{split}
\frac{dy_i}{d\tau}&=\frac{1}{w}(-\frac{k_ia_{i}^2}{2}(\frac{\partial \hat{J}_i(y)}{\partial y_i}+\mu_i)+O(a_{i}^2)+O(a_{i}^3)),\\
\frac{d\mu_i}{d\tau}&=\frac{1}{w}(y_i+O(a_{i})-d_i-\alpha\sum_{j\in\mathcal{N}_i}(\mu_i-\mu_j)-z_i),\\
\frac{dz_i}{d\tau}&=\frac{1}{w}(\alpha\sum_{j\in\mathcal{N}_i}(\mu_i-\mu_j)),\\
\frac{dn_i}{d\tau}&=\frac{1}{w}(-w_{hi}n_i+w_{hi}\hat{J}_i).
\end{split}
\end{align}

Consider the Lyapunov function $V$. It is derived that
\begin{align*}
\begin{split}
\frac{dV}{d\tau}&\leq -(\phi \underline{K}m-\frac{\underline{K}}{2}-\frac{1}{2})\|y-y^*\|^2-\frac{1}{2}\|z-z^*\|^2\\
&\ \ \ -(\alpha \lambda_2-\frac{\phi^2\underline{K}}{2}-2)\|\mu-\mu^*\|^2-\|n-n^*\|^2\\
&\ \ \ +O(\bar{a})+O(\bar{a}^2)+O(\bar{a}^3),
\end{split}
\end{align*}
where $\bar{a}=\max\{a_{1},\ldots,a_{N}\}$.  Since $\underline{K}>\frac{1}{2\phi m-1}$ and $\alpha>\frac{4+\phi^2m}{2\lambda_2}$ with $\phi>\frac{1}{2m}$, there exists a class $\mathcal{K}$ function $\sigma_1$ such that $dV/d\tau\leq -\sigma_1(\|\Theta\|) +O(\bar{a})+O(\bar{a}^2)+O(\bar{a}^3)$, where $\Theta=\col(y-y^*, \mu-\mu^*,z-z^*,n-n^*)$. Recall the definition of the Lyapunov function. There exist positive constants $M_1$ and $M_2$ such that $M_1\|\Theta\|^2\leq V\leq M_2\|\Theta\|^2$. According to the proof of Theorem 4.18 in \cite{Khalil.2002}, take any $r>0$ such that $\boldsymbol{B}_r\subset D$ and choose $\bar{a}$ to be sufficiently small such that $\xi<\sqrt{\frac{M_1}{M_2}}r$, where $\xi=\sigma_1^{-1}(\varphi(\bar{a}))$ with $\varphi(\bar{a})=2(O(\bar{a})+O(\bar{a}^2)+O(\bar{a}^3))$. Let $\rho_1=M_1r^2$ and $\rho_2=M_2\xi^2$ and define $D_{\rho_2}=\{\Theta\in\boldsymbol{B}_r| V\leq \rho_2\}$ and $D_{\rho_1}=\{\Theta\in\boldsymbol{B}_r| V\leq \rho_1\}$. Then, $\boldsymbol{B}_\xi\subset D_{\rho_2}\subset D_{\rho_1}\subset \boldsymbol{B}_r\subset D$. Suppose that $\Delta=\sqrt{\frac{M_1}{M_2}}r$ such that $\|\Theta(\tau_0)\|\leq \Delta$. Then, $\Theta(\tau)$ must enter $D_{\rho_2}$ in a finite time for $\Theta(\tau_0)\in D_{\rho_1}$ and $\tau>\tau_0$. If $\Theta(\tau_0)\in D_{\rho_2}$, a solution starting inside $D_{\rho_2}$ satisfies
\begin{align*}
\|\Theta(\tau)\|\leq \sqrt{\frac{M_2}{M_2}}\sigma_1^{-1}(\varphi(\bar{a})),
\end{align*}
for all $\tau\geq \tau_0$ as $D_{\rho_2}= \{\Theta\in \boldsymbol{B}_r| \|\Theta\|\leq \sqrt{\frac{M_2}{M_1}}\xi\}$.
It is concluded that for the boundary layer system \eqref{bl} and for all $t\geq t_0\geq 0$, there exists a $\mathcal{KL}$ function $\beta$ such that
\begin{align*}
\|\Theta(t)\|\leq \beta(\|\Theta(t_0)\|,(t-t_0))+v,
\end{align*}
where $v=\sqrt{\frac{M_2}{M_1}}\sigma_1^{-1}(\varphi(\bar{a}))+O(1/w)$ by Lemma 3.4, and Theorem 10.4 in \cite{Khalil.2002}.  Via the singular perturbation theory, after the boundary layer system \eqref{bl} arrives at a neighborhood of the GNE whose neighborhood depends on $v$, the state of the reduced system \eqref{rd} exponentially converges to zero. Thus, we have the conclusion in Theorem \ref{theorem1}.
$\hfill\blacksquare$


\begin{remark}
According to the idea of the classical ESC proposed in \cite{Tan.2006}, parameter $w$ should be large to ensure that  dither signals are fast compared with the extremum seeking scheme. Inspired by \cite{Wang.2016}, the amplitude of the dither signals changes adaptively with the extremum estimation error. Once the costs of all players reach the GNE outcome stably, the amplitude of dither signal exponentially converges to zero. Thus, a small parameter $\delta$ is required to guarantee that the dynamics of dither signal amplitude are slower than the extremum seeking scheme. The cut-off frequencies of high-pass and low-pass filters should be selected lower than the frequency of the dither signal.
\end{remark}

\begin{remark}
Different from local results stated in \cite{Frihauf.2012,Zahedi.2019}, Theorem \ref{theorem1} presents a non-local convergence result as the initial state is not required to stay in a small neighborhood of the GNE. Instead, the distance between the initial state and the GNE is bounded. And the bound depends on the maximum amplitude of dither signals at initial time.
\end{remark}

\section{Special Cases}

\subsection{General Quadratic Cost Functions}
Here, the cost function of player $i$ is a general quadratic form, which is given by
\begin{align} \label{gqcf}
J_i(x_i,x_{-i})=\frac{1}{2}\sum_{j=1}^N \sum_{k=1}^N D_{jk}^i x_jx_k+\sum_{j=1}^Nb_j^ix_j+c^i,
\end{align}
for all $i\in\mathcal{I}$, where $x_i\in \mathbb{R}$ is the strategy of player $i$, $D_{jk}^i$, $b_j^i$ and $c^i$ are constants, $D_{ii}^i>0$, and $D_{jk}^i=D_{kj}^i$. Furthermore, strategies of all players are restricted by the equality constraint, that is,
\begin{equation} \label{cc}
C=\{x\in \mathbb{R}^N: \sum_{i=1}^Nx_i=\sum_{i=1}^N d_i\}.
\end{equation}

\begin{assumption} \label{as3}
The matrix $D$ given as follows
\begin{equation}
D=
\begin{bmatrix}
D_{11}^1 & D_{12}^1 & \cdots & D_{1N}^1\\
D_{21}^2 & D_{22}^2 \\
\vdots & & \ddots \\
D_{N1}^N & \cdots & & D_{NN}^N
\end{bmatrix}
\end{equation}
is strictly diagonally dominant, i.e., $\sum_{j\neq i}^N |D_{ij}^i| \leq |D_{ii}^i|, \forall i\in \mathcal{I}$.
\end{assumption}

\begin{lemma}
The noncooperative game with cost functions \eqref{gqcf} and coupled constraint \eqref{cc} admits a GNE $x^*=[x_1^*,\ldots,x_N^*]^T$ if and only if
\begin{equation} \label{gqc}
\begin{aligned}
D_{ii}^ix_i^*+\sum_{j\neq i}D_{ij}^i x_j^*+b_i^i+\bar{\mu}&=0, \forall i \in \mathcal{I},\\
\sum_{i=1}^N x_i^*&= \sum_{i=1}^N d_i,
\end{aligned}
\end{equation}
admits a solution, where $\bar{\mu}\in\mathbb{R}$.
\end{lemma}

Proof: Given the strategies of other players $x_{-i}^*$, player $i$ sloves a constrained minimization problem,
\begin{equation} \label{gqo}
\begin{array}{l}
\min_{x_i\in \mathbb{R}} J_i(x_i,x_{-i}^*)\\
s.t.\ \   x_i+\sum_{j\neq i}^Nx_j^*=\sum_{i=1}^Nd_i.
\end{array}
\end{equation}
Since the cost function $J_i(x,x_{-i}^*)$ is a general quadratic function \eqref{gqcf}, the Lagrange function is constructed by $L(x_i,\mu_i)=J_i(x_i,x_{-i}^*)+\mu_i(x_i+\sum_{j\neq i}^Nx_j^*-\sum_{i=1}^Nd_i)$, where $\mu_i$ is the Lagrange multiplier corresponding to the coupled equality constraint. Due to the convexity of cost function $J_i(x_i,x_{-i}^*)$ with respect to $x_i$, Lagrange function $L(x_i,\mu_i)$ is convex with respect to $x_i$ and linear with respect to $\mu_i$. The solution to \eqref{gqo} is denoted by $(x_i^*,\mu_i^*)$, satisfying the following KKT conditions for all $i\in \mathcal{I}$.
\begin{equation}
\begin{array}{l}
D_{ii}^ix_i^*+\sum_{j\neq i}D_{ij}^i x_j^*+b_i^i+\mu_i^*=0,\\
\sum_{i=1}^Nx_i^*=\sum_{i=1}^N d_i.
\end{array}
\end{equation}
Here, we focus on the variational equilibrium. According to Theorem 3.1 in \cite{Facchinei.2007}, if $\mu_1^*=\ldots =\mu_N^*=\bar{\mu}$ for some $\bar{\mu}\in \mathbb{R}$, the GNE is a variational equilibrium. Conversely, if the GNE is a variational equilibrium, there exists a common $\bar{\mu}\in\mathbb{R}$ such that \eqref{gqc} holds.
$\hfill\blacksquare$

For the general quadratic cost functions \eqref{gqcf} and coupled equality constraint \eqref{cc}, the GNE seeking algorithm is similar to \eqref{sur} with $\hat{J}_i$ replaced by $J_i$ for all $i\in\mathcal{I}$.

\begin{theorem} \label{theorem2}
Let $D$ be a compact set that contains the origin. Suppose that Assumption \ref{as3} holds, the parameters $k$ and $\alpha$ satisfy that $\lambda_{min}(Q)>\frac{1}{2}+\|k^TP\|$ and $\alpha>\frac{2+\|k^TP\|}{\lambda_2}$, respectively, and the frequencies of dither signals satisfy that $w_i\neq w_j$, $w_i\neq w_j+w_k$, and $2w_i\neq w_j$ for any distinct $i,j,k\in\mathcal{I}$. All the players with cost functions \eqref{gqcf} follow the seeking algorithm \eqref{sur} to update their strategies. Then, there exist $w^*>0$ and $\delta^*$ such that for any $w>w^*$, $\delta\in(0,\delta^*)$,
\begin{align*}
\|\Xi(t)\|&\leq \beta(\|\Xi(t_0)\|,(t-t_0))+O(1/w), \ \ \forall t\geq t_0 \geq 0
\end{align*}
with $\Xi(t_0)\in D$, $\|\Xi(t_0)\|\leq\Delta$ and $a_i(t_0)>0, i\in\mathcal{I}$, where $\beta$ is a $\mathcal{KL}$ function, $\Xi(t)=\col(x-x^*,\mu-\mu^*,z-z^*,a-a^*,n-n^*)$, and $\Delta$ depends on $\max\{a_1(t_0),\ldots,a_N(t_0)\}$.
\end{theorem}

Proof: The analysis of the reduced system and the boundary layer system of \eqref{sur} is similar to the proof of Theorem \ref{theorem1}. Substituting the cost function \eqref{gqcf} into the average system \eqref{ave} yields that
\begin{align} \label{rsp}
\begin{split}
\frac{d\hat{x}_i^a}{d\tau}&=\frac{1}{w}(-\frac{k_ia_i^2}{2}(\sum_{j=1}^ND_{ij}^i\hat{x}_j^a+b_i^i)-\frac{k_ia_i^2}{2}\mu_i),\\
\frac{d\mu_i}{d\tau}&=\frac{1}{w}(\hat{x}_i^a-d_i-\alpha\sum_{j\in\mathcal{N}_i}(\mu_i-\mu_j)-z_i),\\
\frac{dz_i}{d\tau}&=\frac{\alpha}{w}\sum_{j\in\mathcal{N}_i}(\mu_i-\mu_j),\\
\frac{dn_i}{d\tau}&=\frac{1}{w}(-w_{hi}n_i+w_{hi}J_i).
\end{split}
\end{align}

Let $\hat{x}^a$ $=$ $[\hat{x}_1^a,\ldots,\hat{x}_N^a]^T$, $\mu$ $=$ $[\mu_1,\ldots,\mu_N]^T$, $z$ $=$ $[z_1,\ldots,$ \\$z_N]^T$, $n=[n_1,\ldots,n_N]^T$, $k$ $=$ $\frac{1}{2}\diag\{k_1a_1^2,\ldots,k_Na_N^2\}$, $b$ $=$ \\ $\frac{1}{2}[k_1a_1b_1^1,\ldots,k_Na_Nb_N^N]^T$, $d$ $=$ $[d_1,\ldots,d_N]^T$, $w_h$ $=$ $\diag\{w_{h1},\ldots,w_{hN}\}$ and $J$ $=$ $[J_1,\ldots,J_N]^T$. The compact form of \eqref{rsp} can be written as follows.
\begin{align} \label{rblc}
\begin{split}
\frac{d\hat{x}^a}{d\tau}&=\frac{1}{w}(A\hat{x}^a-b-k\mu),\\
\frac{d\mu}{d\tau}&=\frac{1}{w}(\hat{x}^a-d-\alpha L\mu-z),\\
\frac{dz}{d\tau}&=\frac{1}{w}\alpha L\mu,\\
\frac{dn}{d\tau}&=\frac{1}{w}(-w_hn+w_hJ),
\end{split}
\end{align}
where
\begin{align*}
A=-\frac{1}{2}
\begin{bmatrix}
k_1a_1D_{11}^1 & k_1a_1D_{12}^1 & \cdots & k_1a_1D_{1N}^1\\
k_2a_2D_{21}^2 & k_2a_2D_{22}^2 \\
\vdots & & \ddots \\
k_Na_ND_{N1}^N & \cdots & & k_Na_ND_{NN}^N
\end{bmatrix}.
\end{align*}
Since $A$ is strictly diagonally dominant, it follows from the Gershgorin Circle Theorem that the eigenvalues of $A$ have negative real parts. Then, there exists a matrix $P=P^T>0$ such that $A^TP+PA=-Q$.

Let $(x^*,\mu^*,z^*,n^*)$ be the equilibrium of system \eqref{rblc}. Denote $\tilde{x}=\hat{x}^a-x^*, \tilde{\mu}=\mu-\mu^*$, $\tilde{z}=z-z^*$, and $\tilde{n}=n-n^*$ to transform the equilibrium to the origin. It yields that
\begin{align} \label{od}
\begin{split}
\frac{d\tilde{x}}{d\tau}&=\frac{1}{w}(A\tilde{x}-k\tilde{\mu}),\\
\frac{d\tilde{\mu}}{d\tau}&=\frac{1}{w}(\tilde{x}-\alpha L\tilde{\mu}-\tilde{z}),\\
\frac{d\tilde{z}}{d\tau}&=\frac{1}{w}\alpha L\tilde{\mu},\\
\frac{d\tilde{n}}{d\tau}&=-\frac{1}{w}w_h\tilde{n}.
\end{split}
\end{align}

Define Lyapunov function $V=w\tilde{x}^TP\tilde{x}+\frac{w}{2}\tilde{\mu}^T\tilde{\mu}+\frac{w}{2}(\tilde{\mu}+\tilde{z})^T(\tilde{\mu}+\tilde{z})+\frac{w}{2}\tilde{n}^Tw_h^{-1}\tilde{n}$. Its derivative along the solutions of \eqref{od} is given by
\allowdisplaybreaks[4]
\begin{align*}
\begin{split}
\dot{V}&=\tilde{x}^T(A^TP+PA)\tilde{x}-2\tilde{\mu}^Tk^TP\tilde{x}+\tilde{\mu}^T\tilde{x}\\
&\ \ \ -\alpha\tilde{\mu}^TL\tilde{\mu}-\tilde{\mu}^T\tilde{z}+\tilde{\mu}^T\tilde{x}-\tilde{\mu}^T\tilde{z}-\|\tilde{z}\|^2-\|\tilde{n}\|^2\\
&\leq -\lambda_{min}(Q)\|\tilde{x}\|^2+\|2k^TP\|\|\tilde{\mu}\|\|\tilde{x}\|-\alpha\lambda_2\|\tilde{\mu}\|^2\\
&\ \ \ -\|\tilde{n}\|^2-\|\tilde{z}\|^2-\|-\sqrt{2}\tilde{\mu}-\frac{\sqrt{2}}{2}\tilde{z}+\frac{\sqrt{2}}{2}\tilde{x}\|^2\\
&\ \ \ +2\|\tilde{\mu}\|^2+\frac{1}{2}\|\tilde{z}\|^2+\frac{1}{2}\|\tilde{x}\|^2\\
&\leq -(\lambda_{min}(Q)-\frac{1}{2}-\|k^TP\|)\|\tilde{x}\|^2\\
&\ \ \ -(\alpha\lambda_2-2-\|k^TP\|)\|\tilde{\mu}\|^2-\frac{1}{2}\|\tilde{z}\|^2-\|\tilde{n}\|^2.
\end{split}
\end{align*}
If $\lambda_{min}(Q)>\frac{1}{2}+\|k^TP\|$ and $\alpha\lambda_2>2+\|k^TP\|$, $\dot{V}<0$.  The average system \eqref{rblc} can asymptotically converge to $(x^*,\mu^*,z^*,n^*)$. Like the analysis of the boundary layer system in Theorem 1, the boundary layer system of \eqref{sur} asymptotically converges to a neighborhood of $(x^*,\mu^*,z^*,n^*)$. By the singular perturbation theory, there exist $w^*$ and $\delta^*$ such that for any $w>w^*$ and $\delta\in (0,\delta^*)$, the system \eqref{sur} asymptotically converges to a small neighborhood of $(x_i^*,\bar{\mu},z_i^*,0,n_i^*)$ for all $i\in \mathcal{I}$.
$\hfill\blacksquare$

\subsection{General Quadratic Games with Stubborn Players}
We assume that player $i\in\{1,\ldots,n\}$ executes \eqref{sur} and player $j\in\{n+1,\ldots, N\}$ is stubborn. Then, player $j$ uses the fixed strategy, i.e., $x_j(t)=\bar{x}_j, \forall j\in \{n+1,\ldots, N\}$. Each player has a quadratic cost function \eqref{gqcf}. In the  presence of stubborn players, the other players' best response $x^{br}=[x_1^{br},\ldots,x_n^{br}]$ is given by
\begin{align}
\begin{split}
\begin{bmatrix}
x_1^{br}\\ \vdots \\ x_n^{br}
\end{bmatrix}
&=-
\begin{bmatrix}
D_{11}^1 & \cdots & D_{1n}^1\\
\vdots & \ddots & \\
D_{n1}^n & \cdots & D_{nn}^n
\end{bmatrix}^{-1}\\
&\ \ \ \times\Big(
\begin{bmatrix}
D_{1,n+1}^1 & \cdots & D_{1N}^1\\
\vdots & \ddots & \\
D_{n,n+1}^n & \cdots & D_{nN}^n
\end{bmatrix}
\begin{bmatrix}
\bar{x}_{n+1}\\
\vdots\\
\bar{x}_{N}
\end{bmatrix}\\
&\ \ \ +
\begin{bmatrix}
b_1^1+\bar{\mu}\\
\vdots\\
b_n^n+\bar{\mu}
\end{bmatrix} \Big).
\end{split}
\end{align}

\begin{theorem} \label{theorem3}
Let $D$ be a compact set that contains the origin. Suppose that Assumption \ref{as3} holds, the parameters $k$ and $\alpha$ satisfy that $\lambda_{min}(Q)>\frac{1}{2}+\|k^TP\|$ and $\alpha>\frac{2+\|k^TP\|}{\lambda_2}$, respectively, and the frequencies of dither signals satisfy $w_i\neq w_j$, $w_i\neq w_j+w_k$, and $2w_i\neq w_j$ for any distinct $i,j,k\in\mathcal{I}$. Player $i$ ($i\in\{1,\ldots,n\}$) with cost function \eqref{gqcf} follows the seeking algorithm \eqref{sur} to update its strategy and player $j$ $(j\in\{n+1,\ldots,N\})$ is stubborn. Then, there exist $w^*$ and $\delta^*$ such that for any $w>w^*$, $\delta\in (0,\delta^*)$,
\begin{align*}
\|\Omega(t)\|\leq \beta(\|\Omega(t_0)\|,(t-t_0))+O(1/w),\ \ \forall t\geq t_0\geq 0
\end{align*}
with $\Omega(t_0)\in D$, $\|\Omega(t_0)\|\leq \Delta$, and $a_i(t_0)>0,i\in\{1,\ldots,n\}$, where $\beta$ is a $\mathcal{KL}$ function, $\Omega(t)=\col(x-x^{br},\mu-\mu^*,z-z^*,a-a^*,n-n^*)$, and $\Delta$ depends on $\max\{a_1(t_0),\ldots,a_n(t_0)\}$.
\end{theorem}

Proof: Similar to the proof of Theorem 2, we can obtain that
\begin{align}
\begin{split}
\frac{d}{d\tau}
\!\begin{bmatrix}
\tilde{x}_1\\
\vdots\\
\tilde{x}_n
\end{bmatrix}
&\!=\!-\frac{1}{2w}\!
\begin{bmatrix}
k_1a_1D_{11}^1 & \cdots & k_1a_1D_{1n}^1\\
\vdots & \ddots & \vdots\\
k_na_nD_{n1}^n & \cdots & k_na_nD_{nn}^n
\end{bmatrix}
\begin{bmatrix}
\tilde{x}_1\\
\vdots\\
\tilde{x}_n
\end{bmatrix}\\
&\ \ \ -\frac{1}{w}k\tilde{\mu},\\
\frac{d\tilde{\mu}}{d\tau}&=\frac{1}{w}(\tilde{x}-\alpha L\tilde{\mu}-\tilde{z}),\\
\frac{d\tilde{z}}{d\tau}&=\frac{1}{w}\alpha L\tilde{\mu},\\
\frac{d\tilde{n}}{d\tau}&=-\frac{1}{w}w_h\tilde{n},
\end{split}
\end{align}
where $\tilde{x}_i=\hat{x}_i^a-x_i^{br}$, $\tilde{x}=[\tilde{x}_1,\ldots, \tilde{x}_n]^T$ $\tilde{\mu}$, $\tilde{z}$ and $\tilde{n}$ are defined in \eqref{od}. So from the proof of Theorem 2, we have that the strategies of $n$ players converge to the best response $x^{br}$.
$\hfill\blacksquare$

\section{Simulations}
In this section, we consider a network of Cournot competition with local and coupled constraints. Four companies produce the same product to meet the market demand. The network of communication among companies is a ring graph. The cost function of company $i$ is $J_i(x_i,x_{-i})=x_i^2+\alpha_ix_i+\beta_i-(p_0-c\sum_{i=1}^4x_i)x_i$, for all $i\in\{1,2,3,4\}$,  where $[\alpha_1,\ldots,\alpha_4]=[2,4,6,8]$, $[\beta_1,\ldots,\beta_4]=[7.5,9,12,15]$, $p_0=5$ and $c=0.04$. The local constraints for the four companies are $1\leq x_1\leq 7$, $2\leq x_2\leq 6$, $2\leq x_3\leq 5$ and $0\leq x_4\leq 4$. The coupled constraint is $x_1+x_2+x_3+x_4=16$. It is calculated that the GNE is $x^*=[5.6,4.5,3.5,2.4]^T$.

\begin{figure}[H]
  \centering
  \includegraphics[width=8cm]{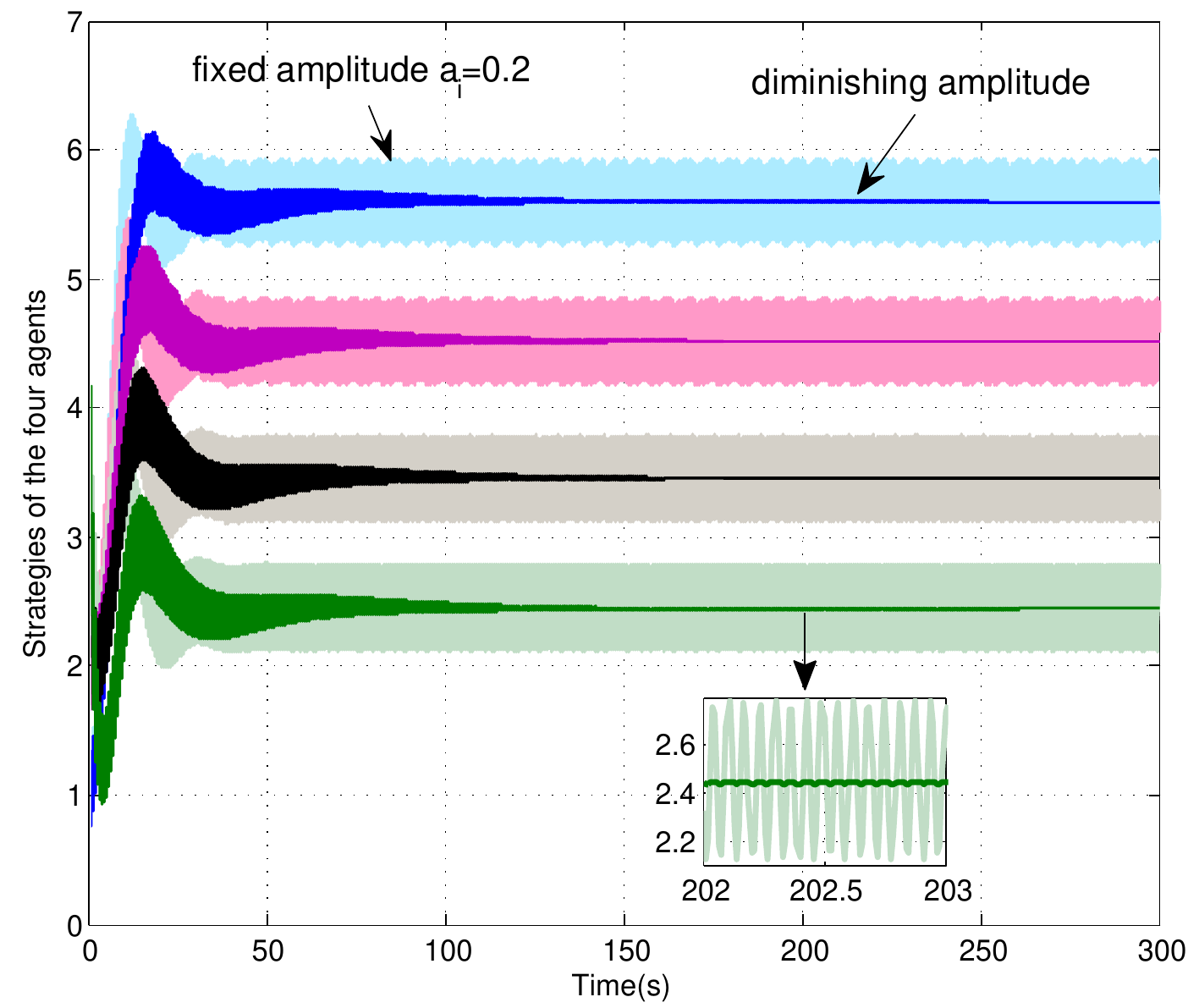}
  \caption{The evolution of the strategies of all agents following \eqref{sur}}
  \label{fig1}
\end{figure}

\begin{figure}[!t]
  \centering
  \includegraphics[width=8cm]{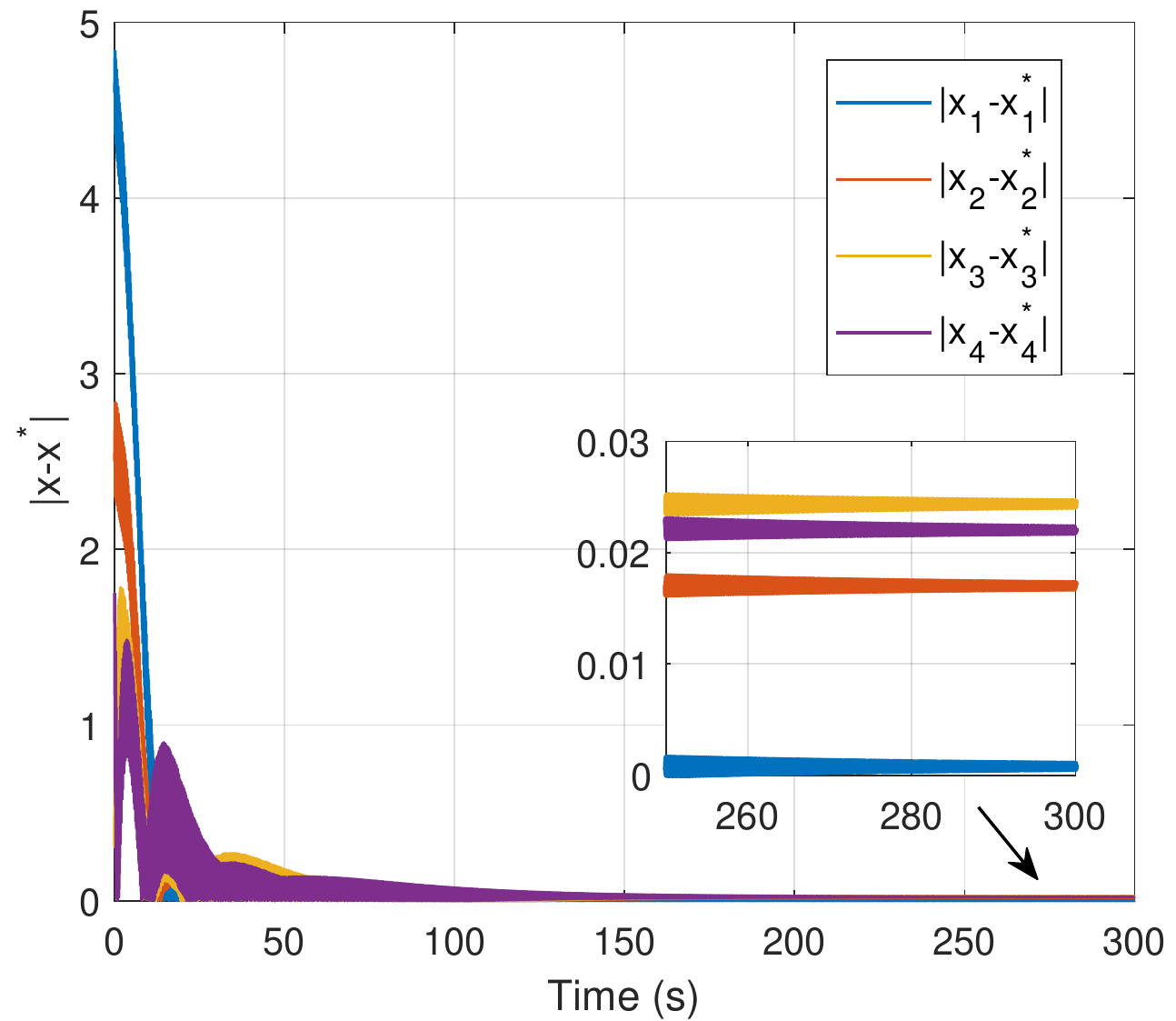}
  \caption{The performance of algorithm \eqref{sur}}
  \label{fig4}
\end{figure}

\begin{figure}[!t]
  \centering
  \includegraphics[width=8cm]{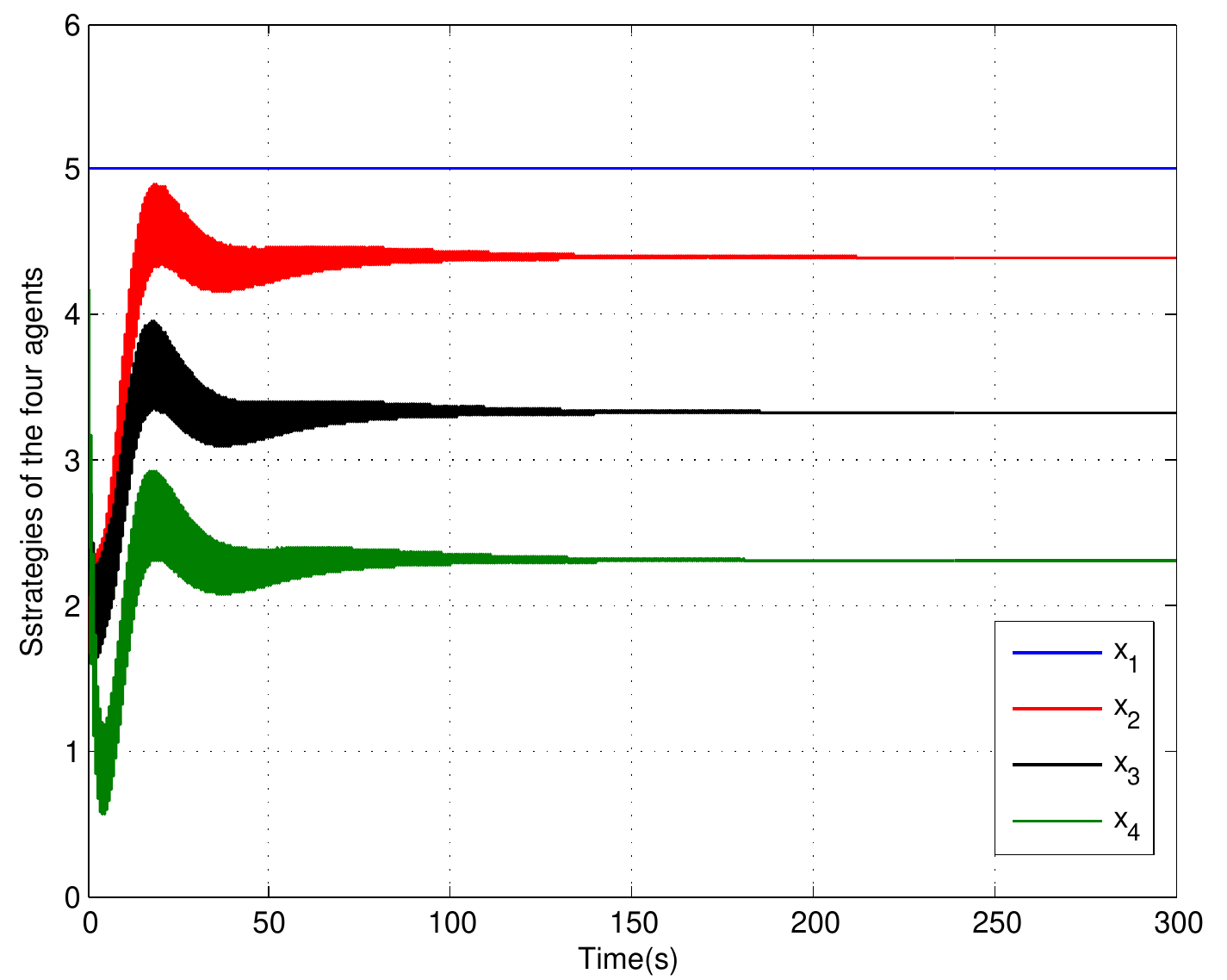}
  \caption{The evolution of the strategies of all agents following \eqref{sur} with one stubborn player}
  \label{fig2}
\end{figure}

Select the parameter $w_{li}=0.5$, $w_{hi}=5$, $b_i=1$ and $k_i=3$ for all $i\in \{1,2,3,4\}$, and $\delta=0.05$. The frequencies of the sinusoidal dither signals for the four agents are selected as $[w_1,w_2,w_3,w_4]=[100,101,103,98]$. The initial conditions are $x(0)=[1,2,3,4]^T$, $a_i(0)=0.2$, $\mu=[0,0,0,0]^T$, $z=[0,0,0,0]^T$ and $n=[0,0,0,0]^T$. The evolution of the strategies of the four players is depicted in Fig.~\ref{fig1}. Compared with the case of the dither signals with the fixed amplitude, the strategies of all players reach to the GNE of the noncooperative game without steady-state oscillations. The performance of the designed algorithm is shown in Fig.~\ref{fig4}.
As a special case, we assume that the company $1$ is a stubborn player whose strategy is $x_1=5$. The best responses of other players are shown in Fig.~\ref{fig2}.

\section{Conclusions}
In this paper, a distributed extremum-seeking based seeking algorithm has been proposed for noncooperative games with coupled constraints to seek GNE. The exact penalty method has been utilized to deal with local constraints to construct auxiliary cost functions whose values can be obtained by measurements of agents' cost functions and local constraints. Thus, the resulting auxiliary problem can be solved by the designed seeking algorithm without any knowledge of the explicit expressions of cost functions and local constraints. The coupled constraints that limit the strategies of all agents have been handled by the coordination of Lagrange multipliers. Via singular perturbation theory, averaging analysis and Lyapunov stability theory, the non-local convergence of all agents' strategies has been analyzed under certain conditions. Additionally, it is interesting to employ the projection operator to deal with local box constraints which may be nonconvex.


%



\ifCLASSOPTIONcaptionsoff
  \newpage
\fi



%


\end{document}